\documentclass[10pt]{amsart}
\usepackage{latexsym, amsmath}

 \renewcommand{\epsilon}{\varepsilon}
 \newcommand{\newsection}[1]
 {\subsection{#1}\setcounter{theorem}{0} \setcounter{equation}{0}
 \par\noindent}

 \newtheorem{theorem}{Theorem}
 
 \newtheorem{lemma}[theorem]{Lemma}
 \newtheorem{corr}[theorem]{Corollary}
 
 \newtheorem{proposition}[theorem]{Proposition}
 \newtheorem{deff}[theorem]{Definition}
 
 \newcommand{\bth}{\begin{theorem}}
 \newcommand{\ble}{\begin{lemma}}
 \newcommand{\bcor}{\begin{corr}}
 \newcommand{\bdeff}{\begin{deff}}
 \newcommand{\bprop}{\begin{proposition}}
 \newcommand{\eth}{\end{theorem}}
 \newcommand{\ele}{\end{lemma}}
 \newcommand{\ecor}{\end{corr}}
 \newcommand{\edeff}{\end{deff}}
 
 \newcommand{\eprop}{\end{proposition}}
  \newcommand{\cd}{\, \cdot\, }

 \renewcommand{\Pi}{\varPi}

 \renewcommand{\epsilon}{\varepsilon}

 \newcommand{\tidle}{\tilde}

\begin{document}

 \title{Eigenfunction and Bochner Riesz estimates on manifolds with boundary}

\thanks{The author was supported in part by the NSF}
\author{Christopher D. Sogge}
\address{Department of Mathematics, The Johns Hopkins University,
Baltimore, MD 21218}

 \maketitle

 \newsection{Introduction}

 The purpose of this paper is to give a simple proof of sharp
 $L^\infty$ estimates for the eigenfunctions of the Dirichlet
 Laplacian on smooth compact Riemannian manifolds $(M,g)$ of dimension $n\ge 2$ with
 boundary $\partial M$ and then to use these estimates to prove new estimates for
 Bochner Riesz means in this setting.  Thus, we shall consider the Dirichlet
 eigenvalue problem
 \begin{equation}\label{0.0}
 (\Delta +\lambda^2)u(x)=0, x \in M, \quad u(x)=0, \, \,
 x\in \partial M,
  \end{equation}
 with $\Delta=\Delta_g$ being the Laplace-Beltrami operator
 associated to the Riemannian metric $g$.
 Recall that the spectrum of $-\Delta$ is discrete and tends to
 infinity.  Let $0<\lambda_1^2\le \lambda^2_2\le \lambda_3^2\le
 \dots$ denote the eigenvalues, so that $\{\lambda_j\}$ is the
 spectrum of the first order operator
 $\mathcal{P}=\sqrt{-\Delta}$.  Let $\{e_j(x)\}$ be an associated real
 orthonormal basis, and let
 $$e_j(f)(x)=e_j(x) \int_M f(y) e_j(y) \, dy,$$
 be the projection onto the $j$-th eigenspace.  Here and in what
 follows $dy$ denotes the volume element associated with the
 metric $g$.

Grieser \cite{G2} recently proved the $L^\infty$ estimates
 \begin{equation}\label{0.1}
 \|e_j(f)\|_\infty \le C\lambda_j^{(n-1)/2}\|f\|_2,
 \end{equation}
 which are sharp for instance when $M$ is the upper hemisphere of
 $S^n$ with the standard metric.  One of our main results is a
 slight strengthening of this.  We shall consider the unit band
 spectral projection operators,
 \begin{equation}\label{0.2}
 \chi_\lambda f = \sum_{|\lambda_j-\lambda|\le 1} e_j(f),
 \end{equation}
 and show that these enjoy the same bounds:

 \begin{theorem}\label{theorem1}  Fix a compact Riemannian
 manifold $(M,g)$ with boundary of dimension $n\ge 2$.  Then there
 is a uniform constant $C$ so that
 \begin{equation}\label{0.3}
 \|\chi_\lambda f\|_\infty \le C \lambda^{(n-1)/2}\|f\|_2, \quad
 \lambda \ge 1.
 \end{equation}
 \end{theorem}

 In the case of manifolds without boundary, this and more general
 estimates of the form
\begin{equation}\label{0.4}\|\chi_\lambda f\|_p \le C
 \lambda^{\sigma(p)}\|f\|_p, \quad \lambda \ge 1, \, p\ge 2
 \end{equation}
 where
 \begin{equation}\label{0.5}
 \sigma(p)=\max\left\{\,  \frac{n-1}2-\frac{n}p, \,
 \frac{n-1}2\left(\frac12-\frac1p\right)  \, \right\}
 \end{equation}
 were proved in \cite{S1}.  These estimates cannot be improved
 since one can show that the operator norms satisfy
 $\limsup_{\lambda \to\infty}\lambda^{-\sigma(p)}\|\chi_\lambda\|_{L2\to
 L^p}>0$  (see \cite{S3}).

 The special case of \eqref{0.4} where $p=\infty$
 seems to have been first stated in \cite{S1}, but it can be proved using
 much older estimates of Avakumovic \cite{A1}, \cite{A2} and Levitan \cite{L}
 that arise in the proof of the sharp Weyl formula for Riemannian
 manifolds without boundary.  After that, H\"ormander \cite{H1}
 proved the sharp Weyl formula for general self-adjoint elliptic
 operators on manifolds without boundary.  Recently, in the case
 of manifolds without boundary, the author and Zelditch \cite{SZ} proved
 estimates that imply that for generic metrics on any manifold one
 has the bounds $\|e_j\|_\infty=o(\lambda_j^{(n-1)/2})$ for
 $L^2$-normalized eigenfunctions.  The corresponding result for
 manifolds with boundary is not known.

 In the case of manifolds with
 boundary, the only known results for the unit band spectral
 projection operators were due to D. Grieser \cite{G1} and H. Smith and the author
 \cite{SS}, who showed that the bounds \eqref{0.4} hold under the
 assumption that the manifold has geodesically concave boundary.
The two-dimensional case was handled in \cite{G1}, and the higher
dimensional in \cite{SS}.

 In the other direction, Grieser in his thesis \cite{G1}, showed
 that \eqref{0.4} cannot hold if the boundary of $M$ has a point
 that is geodesically convex.  In this case, for instance, when
 $n=2$, he constructed a counterexample showing that the bounds in
 in \eqref{0.4} can only hold for $p\ge 8$.  Showing that they are
 valid for this range of exponents appears to be very difficult.
 The reason for the difference in the smaller range of exponents
 for this case is related to the existence of Rayleigh whispering
 galleries.  Specifically, one can construct functions with
 spectrum in $\lambda$-unit bands that are concentrated in a
 $\lambda^{-2/3}$ neighborhood of the boundary, while in the
 boundaryless case, the counterexample showing that \eqref{0.4} is
 sharp involves showing that there are functions of this type
 concentrating in $\lambda^{-1/2}$ neighborhoods of geodesics.

We shall follow an idea of Grieser \cite{G2} to prove our
generalization \eqref{0.3} of his result \eqref{0.1}.  Grieser
first showed that one always has the uniform bounds $|e_j(x)|\le
C\lambda^{(n-1)/2}_j$ when the distance from $x$ to $\partial M$
is bounded from below by $\lambda^{-1}_j$.  This first step was
achieved by adapting the proof of estimates for the local Weyl law
which are due to Seeley \cite{SE}, Pham The Lai \cite{P} and
H\"ormander \cite{H2}.  He then used these bounds and a form of
the maximum principle (\cite{PW}, Theorem 10, p. 73) for solutions
of \eqref{0.0} to obtain the bounds in the $\lambda^{-1}_j$
neighborhood of the boundary.

Our proof of \eqref{0.3} will be to first see that the
aforementioned local Weyl estimate of Seeley, Pham The Lai, and
H\"ormander immediately gives the stronger estimate
\begin{equation}\label{0.6}
|\chi_\lambda f(x)|\le C\lambda^{(n-1)/2}\|f\|_2, \quad \lambda\ge
1,\end{equation} for all $x$ when $n=2$ or when for $n\ge 3$ the
distance to the boundary is bounded below by
$\lambda^{-(n-1)/(n-2)}$.  If $n=3$, one can use this fact and a
simple argument involving Sobolev's theorem (cf. Theorem 17.5.3 in
\cite{H2}) to show that \eqref{0.6} must also hold in the missing
piece where $\text{dist}(x,\partial M)\le \lambda^{-2}$ . For
$n\ge 4$; however, one must use a maximum principle argument as
introduced by Greiser for these problems. We shall not directly
use the form of the maximum principle employed by Grieser, but
rather see that its proof can be used to show that the uniform
bounds \eqref{0.6} must also hold in a $\lambda^{-1}$ neighborhood
of the boundary, which would finish the proof of \eqref{0.3}.

Our other main result will be some new estimates for Bochner Riesz
means of eigenfunctions.  Recall the Bochner Riesz means of index
$\delta \ge 0$ are defined by
\begin{equation}\label{0.7}
S^\delta_\lambda f = \sum_{\lambda_j\le \lambda}
\left(1-\frac{\lambda_j^2}{\lambda^2}\right)^\delta e_j(f).
\end{equation}
It is known that a necessary condition for these operators to be
uniformly bounded on $L^p$ for a given $1\le p\le \infty$, $p\ne
2$, is that
\begin{equation}\label{0.8}
\delta(p)=\max\left\{n|1/2-1/p|-1/2, 0\right\}.
\end{equation}
Note that when $p\ge 2(n+1)/(n-1)$, one has $\delta(p)=\sigma(p)$,
where $\sigma(p)$ is the exponent in \eqref{0.5}.  Using this
fact, the author used the boundarlyless estimates \eqref{0.4} to
prove the first sharp estimates for Bochner Riesz means on compact
Riemannian manifolds in \cite{S2}.  Specifically, it was shown
that for a given $p\in [1,2(n+1)/(n+3)] \cup
[2(n+1)/(n-1),\infty]$ one has the uniform bounds
$$\|S_\lambda^\delta f\|_p\le C\|f\|_p$$
in this setting, provided that $\delta>\delta(p)$.  Earlier,
weaker results were due to many people, including H\"ormander
\cite{H3}.

Since we only know that the desired bounds for eigenfunctions
\eqref{0.4} hold for $p=\infty$, we can only at this stage prove
the sharp estimates for Bochner Riesz means when $p=1$ or
$p=\infty$:

\begin{theorem}\label{theorem2}  Fix a smooth compact Riemannian
manifold with boundary $(M,g)$ of dimension $n\ge 2$.  Then if
$\delta>(n-1)/2$ one has the uniform bounds
\begin{equation}\label{0.9}
\|S^\delta_\lambda f\|_p\le C\|f\|_p,
\end{equation}
for every $1\le p\le \infty$.
\end{theorem}

By interpolating with the trivial estimate for $p=2$, and using
duality one gets the bounds \eqref{0.9} from the special case
where $p=1$.  However, the bounds for $1<p<\infty$ certainly are
not sharp.

We shall adapt the argument from \cite{S2} to show that
\eqref{0.4} implies \eqref{0.9}.  In \cite{S2} the Tauberian
arguments behind the proof of the sharp Weyl formula were adapted
to show that one could write $S^\delta_\lambda = \tilde
S^\delta_\lambda + R^\delta_\lambda$, where the "remainder" term
$R^\delta_\lambda$ could be controlled by \eqref{0.4}, while the
other piece, $\tilde S^\delta_\lambda$, could be estimated by
computing its kernel explicitly via the Hadamard parametrix and
then estimating the resulting integral operator using
straightforward adaptations of the arguments for the Euclidean
case.  In the setting of manifolds with boundary, this approach
does not seem to work since the known parametrices for the wave
equation do not seem strong enough unless one assumes that the
boundary is geodesically concave.  Here, we shall get around this
fact by simplifying the earlier arguments and show that estimates
for Bochner Riesz operators just follow from \eqref{0.4} and the
finite propagation speed of solutions of the Dirichlet wave
equation.

In what follows we shall use the convention that $C$ will denote a
constant that is not necessarily the same at each occurrence.

It is a pleasure to thank S. Zelditch for a number helpful
conversations.  I am also grateful to M. Taylor and X. Xu for
helpful comments regarding an early draft of this paper.

\newsection{$L^\infty$ estimates for unit band spectral projection
operators}

In this section we shall prove Theorem \ref{theorem1}.  Thus, we
need to see that one has the uniform bounds
\begin{equation}\label{2.1}|\chi_\lambda f(x)|\le
C\lambda^{(n-1)/2}\|f\|_2, \quad \lambda\ge 1.
\end{equation}
Note that
$$\chi_\lambda f(x)= \int_M \sum_{|\lambda_j-\lambda|\le 1}
e_j(x)e_j(y) f(y) \, dy,$$ therefore, by the converse to
Schwartz's inequality and orthogonality, one has the bounds
\eqref{2.1} at a given point $x$ if and only if
\begin{equation}\label{2.2}
\sum_{|\lambda-\lambda_j|\le 1}(e_j(x))^2 \le C^2 \lambda^{n-1}.
\end{equation}
Because of this, \eqref{2.1} would be a consequence of the
following two results.

\begin{proposition}\label{prop1}  Fix $(M,g)$ as above.  Then,
given $\varepsilon>0$, there is a uniform constant $C$ so that for
$\lambda\ge 1$
\begin{equation}\label{2.3}
\sum_{|\lambda-\lambda_j|\le 1}(e_j(x))^2 \le C_\varepsilon
\lambda^{n-1},
\end{equation}
%for all $x$ if $n=2$, and
for $x$ satisfying
\begin{equation}\label{2.4}
d(x)\ge \varepsilon \lambda^{-1}
%\lambda^{-(n-1)/(n-2)},
\end{equation}
if $n\ge 3$ where $d(x)$ denotes the Riemannian distance to
$\partial M$.
\end{proposition}

\begin{proposition}\label{prop2}  If $(M,g)$ is as above then
for large $\lambda$ we have
\begin{equation}\label{2.5}
\max_{\{x: \, d(x)\le\frac12
(\lambda+1)^{-1}\}}\sum_{|\lambda-\lambda_j|\le 1}(e_j(x))^2 \le 4
\max_{\{x: \, d(x)=\frac12
(\lambda+1)^{-1}\}}\sum_{|\lambda-\lambda_j|\le 1}(e_j(x))^2.
\end{equation}
\end{proposition}

\noindent{\bf Proof of Proposition \ref{prop1}:}  We shall see
that \eqref{2.3} is an immediate consequence of Theorem 17.5.10 in
H\"ormander \cite{H2}, which in turn is based on earlier work of
Seeley \cite{SE} and Pham The Lai \cite{P}.  To state this result,
we let
\begin{equation}\label{2.6}
e(x,\lambda)=(2\pi)^{-n} \int_{\{\xi \in \mathbb{R}^n: \,
%\sum_{i,j=1}^n g^{ij}(x)\xi_i\xi_j
|\xi| \le \lambda\}} \bigl(1-e^{i 2d(x)\xi_n}\bigr)\, d\xi.
\end{equation}
%where $g^{ij}(x)d\xi_id\xi_j$ denotes the cometric, which is dual
%to the Riemannian metric $g_{ij}(x)dx^idx^j$.
If we assume that
local coordinates have been chosen so that the Riemannian volume
form is $dx_1\dots dx_n$, then the result just quoted says that
there is a uniform constant $C$ so that for $\lambda \ge 1$
\begin{equation}\label{2.7}
|\sum_{\lambda_j\le \lambda}(e_j(x))^2-e(x,\lambda)|\le C\lambda
(\lambda+d(x)^{-1})^{n-2}.
\end{equation}
Since $\lambda(\lambda+d(x)^{-1})^{n-2}=O(\lambda^{n-1})$ for all
$x$ satisfying \eqref{2.4}
%if $n=2$ and for $x$ satisfying \eqref{2.5} if $n\ge3$,
this
yields Proposition \ref{prop1} since
$$|e(x,\lambda+1)-e(x,\lambda-1)|\le C\lambda^{n-1}.$$

Although, we do not need to use it here, the proof of \eqref{2.7}
actually gives the bounds \eqref{2.3} when $d(x)\ge
\lambda^{-\frac{n-1}{n-2}}$ for $n\ge 3$ and for all $x$ when
$n=2$.  This stronger fact just follows from the estimate
(17.5.20) in \cite{H2}.

\bigskip

\noindent{\bf Proof of Proposition \ref{prop2}:}  It is convenient
to use geodesic coordinates with respect to the boundary.
Specifically, we shall use the fact that we can find a small
constant $c>0$ so that the map $(x',x_n)\in \partial M \times
[0,c)\to M$, sending $(x',x_n)$ to the endpoint, $x$, of the
geodesic of length $x_n$ which starts at $x'\in \partial M$ and is
perpendicular to $\partial M$ is a local diffeomorphism.  Note
then that $d(x)=x_n$.  Under this identification the metric takes
the form
$$\sum_{i,j=1}^n g_{ij}(x)dx^idx^j=
dx^2_n +\sum_{i,j=1}^{n-1} g'_{ij}(x',x_n) dx^idx^j,$$ where
$g'_{ij}(\, \cdot\, ,x_n)$ is a Riemannian metric on $\partial M$
depending smoothly on $x_n\in [0,c)$.  Consequently, in this
neighborhood of the boundary, the Laplacian can be written as
$$\Delta=\sum_{i,j=1}^n
g^{ij}(x)\frac{\partial^2}{\partial x_i\partial x_j} +\sum_{i=1}^n
b_i(x)\frac{\partial}{\partial x_i},$$ using local coordinates for
$\partial M$, where $g^{ij}$  the matrix with entries
$(g^{ij})_{1\le i,j\le n-1}=(g'_{ij})^{-1}$ and $g^{nn}=1$, and
$g^{nk}=g^{kn}=0$, $k\ne n$.  Also the $b_j(x)$ are $C^\infty$ and
real valued.

In what follows we shall assume that $\lambda$ is large enough so
that $\lambda \ge 2/c$. Assume further that $\text{spec
}\sqrt{-\Delta}\cap [\lambda-1,\lambda+1]\ne \emptyset$, and
consider the function
$$H(x)=\frac1{(w(x))^2}\sum_{\lambda_j\in [\lambda-1,\lambda+1]}(e_j(x))^2,$$
where
$$w(x)=1-(\lambda+1)^2 x_n^2.$$

Suppose that in the strip $\{x\in M: \, 0\le x_n\le\frac12
(\lambda+1)^{-1}\}$ the function $H(x)$ has a maximum at an
interior point $x=x_0$. Then
$$v(x)=\frac1{w(x)}\sum_{\lambda_j\in
[\lambda,\lambda+1]}\frac{e_j(x_0)}{w(x_0)}e_j(x)$$ must have a
positive maximum at $x=x_0$.  For because of our assumptions on
the spectrum we then have $v(x_0)=\sum_{\lambda_j\in
[\lambda,\lambda+1]}(e_j(x_0)/w(x_0))^2>0$, while at other points
in the strip
\begin{multline*}|v(x)|\le
\frac1{w(x)}\bigl(\sum_{\lambda_j\in [\lambda,
\lambda+1]}(e_j(x))^2\bigr)^{1/2} \,
\frac1{w(x_0)}\bigl(\sum_{\lambda_j\in [\lambda,
\lambda+1]}(e_j(x_0))^2\bigr)^{1/2}
\\
= (H(x))^{1/2}(H(x_0))^{1/2}\le H(x_0)=v(x_0).
\end{multline*}

Note that in the strip $\{x\in M: \, 0\le x_n\le\frac12
(\lambda+1)^{-1}\}$ we have
\begin{multline*}(\Delta
+\lambda^2_j)w=-2(\lambda+1)^2-2b_n(x)x_n(\lambda+1)^2+\lambda_j^2(1-(\lambda+1)^2x^2_n)
\\
\le-(\lambda+1)^2/2, \quad \lambda_j\le \lambda+1,
\end{multline*}
assuming that $\lambda$ is large enough so that $|2b_1(x)x_n|\le
1/2$.  Also, in this strip we have that $\frac12\le w(x)\le 1$.

Let us set
$$v_j(x)=\frac{e_j(x)}{w(x)}\frac{e_j(x_0)}{w(x_0)},$$
so that $v(x)=\sum_{\lambda_j\in [\lambda-1,\lambda+1]}v_j(x)$. We
also set
$$u_j(x)=\frac{e_j(x_0)}{w(x_0)}e_j(x),$$
and note that $(\Delta+\lambda^2_j)u_j(x)=0$.

A computation (see p. 72, \cite{PW}) shows that for a given $j$ we
have
\begin{multline*}
0=\frac1w (\Delta+\lambda_j^2)u_j
\\
= \sum_{k,l=1}^n g^{kl}(x)\frac{\partial^2v_j}{\partial
x_k\partial x_l} +\sum_{k=1}^n\bigl(\frac2w
\sum_{l=1}^ng^{kl}(x)\frac{\partial w}{\partial x_l}+b_k\bigr)
\frac{\partial v_j}{\partial x_k} +\frac{v_j}{w}(\Delta
+\lambda_j^2)w.
\end{multline*}
Therefore, if we sum over $\lambda_j\in [\lambda-1,\lambda+1]$, we
get
$$\sum_{k,l=1}^n g^{kl}(x)\frac{\partial^2v}{\partial
x_k\partial x_l} +\sum_{k=1}^n\bigl(\frac2w
\sum_{l=1}^ng^{kl}(x)\frac{\partial w}{\partial x_l}+b_k\bigr)
\frac{\partial v}{\partial x_k} =-\sum_{\lambda_j\in
[\lambda-1,\lambda+1]} \frac{v_j}{w}(\Delta +\lambda_j^2)w.$$ In
particular, at $x=x_0$, we have
\begin{multline*}\sum_{k,l=1}^n
g^{kl}(x_0)\frac{\partial^2v(x_0)}{\partial x_k\partial x_l}
+\sum_{k=1}^n\bigl(\frac2w \sum_{l=1}^ng^{kl}(x_0)\frac{\partial
w}{\partial x_l}+b_k\bigr) \frac{\partial v(x_0)}{\partial x_k}
\\
=\frac{-1}{w(x_0)}\sum_{\lambda_j\in [\lambda-1,\lambda+1]}
\left(\frac{e_j(x_0)}{w(x_0)}\right)^2(\Delta
+\lambda_j^2)w(x_0)>0.
\end{multline*}

But this is impossible since $v$ has a positive maximum at $x_0$,
which implies that $\partial v(x_0)/\partial x_k=0$ for every $k$,
and $\sum_{k,l=1}^n g^{kl}(x_0)\frac{\partial^2v(x_0)}{\partial
x_k\partial x_l}\le 0$. Thus, we conclude that the function $H(x)$
cannot have a maximum at an interior point of the strip, $\{x: \,
0\le x_n \le \frac12(\lambda+1)^{-1}\}$. Because of this, the
Dirichlet conditions, and our lower bound for $w$, we get that
$$\sup_{\{x: \, 0\le x_n\le \frac12(\lambda+1)^{-1}\}}\sum_{\lambda_j\in
[\lambda-1,\lambda+1]}(e_j(x))^2 \le 4\sup_{\{x: \,  x_n=\frac12
(\lambda+1)^{-1}\}}\sum_{\lambda_j\in
[\lambda-1,\lambda+1]}(e_j(x))^2,$$ as desired, which completes
the proof of Proposition \ref{prop2}.

\newsection{Estimates for Bochner Riesz means}

In this section we shall see how favorable estimates for the unit
band spectral projection operators imply sharp estimates for
Bochner Riesz means.  Specifically, we shall prove the following
result which implies Theorem \ref{theorem2}.

\begin{proposition}\label{prop3.1}  Assume that for a given $1\le
p<2$ one has the uniform bounds
\begin{equation}\label{3.1}
 \|\chi_\lambda f\|_2 \le
C\lambda^{\delta(p)}\|f\|_p, \quad \lambda\ge 1,
\end{equation}
where $\delta(p)$ is as in \eqref{0.8}.  Then for a given
$\delta>\delta(p)$ there is a uniform constant $C$ so that
\begin{equation}\label{3.2}
\|S^\delta_\lambda f\|_p\le C_\delta \|f\|_p.
\end{equation}
\end{proposition}

This implies Theorem \ref{theorem2}, since, by duality, (1.4)
implies that \eqref{3.1} must hold when $p=1$.  This implies that
if $\delta>\delta(1)=(n-1)/2$ the $S^\delta_\lambda$ are uniformly
bounded on $L^1$, which implies the same for all $L^p$, $1\le p\le
\infty$ by duality and interpolation.

To prove Proposition \ref{prop3.1}, we shall require the following
straightforward consequences of its hypotheses.

\begin{lemma}\label{lemma3.2}  Suppose that \eqref{3.1} holds.
Suppose also that $\rho\in C(\mathbb{R})$ satisfies
$|\rho(\tau)|\le C_N(1+|\tau|)^{-N}$ for some $N\ge \delta(p)+1$.
Assume also that $1\le 2^k\le \lambda$. Then there is a uniform
constant $C$ so that
\begin{equation}\label{3.3}
\|\rho(2^{-k}(\lambda-\mathcal{P}))f\|_2+\|\rho(2^{-k}(\lambda+\mathcal{P}))f\|_2\le
C2^{k/2}\lambda^{\delta(p)}\|f\|_p, \quad \lambda\ge 1,
\end{equation}
where the constant only depends on $C_N$ and the
constant in \eqref{3.1}.
\end{lemma}

Here we are of course using the notation that $\rho(\mathcal{P})f
= \sum_{j} \rho(\lambda_j)e_j(f)$.

\noindent{\bf Proof of Lemma \ref{lemma3.2}:}  If we just use
orthogonality, \eqref{3.1}, and our assumptions on $\rho$ we find
that
%\begin{align*}
%\|&\rho(2^{-k}(\lambda-\mathcal{P}))f\|_2^2+\|\rho(2^{-k}(\lambda+\mathcal{P}))f\|_2^2
%\\
%&\le C\sum_{j=0}^\infty \sup_{\lambda_l\in
%[j,j+1]}|\rho(2^{-k}(\lambda-\lambda_l))|^2\|\chi_j f\|_2^2 +
%C\sum_{j=0}^\infty \sup_{\lambda_l\in
%[j,j+1]}|\rho(2^{-k}(\lambda-\lambda_l))|^2\|\chi_j f\|_2^2
%\\
%&\le C\sum_{j=0}^\infty (1+2^{-k}|\lambda-j|)^{-N}
%(1+j)^{2\delta(p)}\|f\|_p^2 +C\sum_{j=0}^\infty
%(1+2^{-k}|\lambda+j|)^{-N} (1+j)^{2\delta(p)}\|f\|_p^2.
%\end{align*}
\begin{align*}
\|&\rho(2^{-k}(\lambda-\mathcal{P}))f\|_2^2+\|\rho(2^{-k}(\lambda+\mathcal{P}))f\|_2^2
\\
&\le C\sum_{j=0}^\infty \Bigl(\sup_{\lambda_l\in
[j,j+1]}|\rho(2^{-k}(\lambda-\lambda_l))|^2 +
 \sup_{\lambda_l\in
[j,j+1]}|\rho(2^{-k}(\lambda+\lambda_l))|^2\Bigr)\|\chi_j f\|_2^2
\\
&\le C\sum_{j=0}^\infty \Bigl((1+2^{-k}|\lambda-j|)^{-N}
(1+j)^{2\delta(p)} + (1+2^{-k}|\lambda+j|)^{-N}
(1+j)^{2\delta(p)}\Bigr)\|f\|_p^2.
\end{align*}
The first term in the left dominates the second term. Since
$N>\delta(p)+1$ and $2^{-k}\lambda\ge 1$, by comparing the sums to
the corresponding integrals, one sees that both terms on the right
can be dominated by the square of the right hand side of
\eqref{3.3}, which finishes the proof. \qed

\medskip

We now have the main tools needed to prove Proposition
\ref{prop3.1}.  To be able to rewrite the operators
$S^\delta_\lambda$ in a way that will allow us to use the above
estimates we need to relate the operator to the wave equation. For
this purpose, we need to compute the Fourier transform of the
symbol $\tau \to (1-\tau^2/\lambda^2)^\delta_+$ of the Bochner
Riesz means. We shall use the Bessel function formula,
$$
\int_{-1}^1 e^{i\tau t} (1-\tau^2)^\delta \, d\tau
%=2\int_0^1
%\cos(\tau t) (1-\tau^2)^\delta \, d\tau
=\sqrt{\pi}\Gamma(1+\delta)\left(\frac{t}2\right)^{-\delta-\frac12}J_{\delta+\frac12}(t),
\quad t>0 .
$$
Recall that as $r\to \infty$, we have the following asymptotics
for Bessel functions of order $m$
$$J_m(r)=\sum_{\pm} \alpha_m^\pm(r) e^{\pm ir},$$
where
$$|\partial^j \alpha_m^\pm(r)|\le C_j r^{-j-1/2}, \quad r\ge 1, \, \, j=0,1,2,\dots.$$

Therefore, since this Fourier transform is an even function, we
can write
\begin{align}\label{3.4}
S^\delta_\lambda f&=\sum_{\lambda_j\le
\lambda}(1-\lambda^2_j/\lambda^2)^\delta e_j(f)
\\
&=\frac1{2\pi}\int_{-\infty}^\infty \lambda
\sqrt{\pi}\Gamma(1+\delta)\left(\frac{\lambda
|t|}2\right)^{-\delta-\frac12}J_{\delta+\frac12}(\lambda |t|)
\sum_{j} \cos t\lambda_j \, e_j(f) \, dt \notag
\\ &=\frac1{2\pi}\int_{-\infty}^\infty \sum_{\pm}\lambda
m_\delta^\pm(\lambda t) e^{\pm i \lambda t}
 \cos t\mathcal{P}
f\, dt, \notag
\end{align}
where, for every $j$, $|(1+|t|)^j\partial_t^j m_\delta^\pm(t)|\le
C_\delta (1+|t|)^{-1-\delta}$, and hence
\begin{equation}\label{3.5}
\lambda|(1+|t|)^j\partial^j_t m_\delta^\pm(\lambda t)|\le
\begin{cases}C_\delta \lambda, \quad |t|\le \lambda^{-1}
\\
C_\delta \lambda^{-\delta}|t|^{-1-\delta}, \quad |t|\ge
\lambda^{-1}.
\end{cases}
\end{equation}

Here, $\mathcal{P}=\sqrt{-\Delta}$, and
$$u(t,x)=\cos t\mathcal{P}f(x)=\sum_{j=1}^\infty \cos t\lambda_j e_j(f)(x),$$
is the cosine transform of $f$.  Thus, it is the solution of the
wave equation
$$(\partial_t^2-\Delta_g)u=0, \quad
u(0,\cd)=f, \, \, \partial_t u(0,\cd)=0.$$ We shall use the finite
propagation speed for solutions to the wave equation.
Specifically, if $f$ is supported in a geodesic ball $B(x_0,R)$
centered at $x_0$ with radius $R$, then $x\to\cos t\mathcal{P}f$
vanishes outside of $B(x_0,2R)$ if $0\le t\le R$.

We shall now proceed to break up the operators $S^\delta_\lambda$
into a sum of pieces that we can estimate using a combination of
\eqref{3.3} and H\"older's inequality, along with a ``remainder
term". This is related to an argument of Fefferman \cite{F} for
the Euclidean case, and also an argument of the author \cite{S2}
for the case of Riemannian manifolds without boundary.  The latter
argument also relied on the small time parametrix for the wave
equation, which is impossible to use in this setting.  Instead we
use a simpler argument that only uses the finite propagation speed
of the wave equation.

Let us first deal with the remainder term.  We fix
an even
function
 $b\in C^\infty(\mathbb{R})$ satisfying $b(t)=0$, $|t|<1$
and $b(t)=1$, $|t|>2$, and then set
$$R^\delta_\lambda f
=\frac1{2\pi}\int_{-\infty}^\infty \sum_{\pm}\lambda
m_\delta^\pm(\lambda t)b(t) e^{\pm i \lambda t}
 \cos t\mathcal{P}
f\, dt.$$ If $\rho_\lambda$ denotes the inverse Fourier transform
$t\to \frac12m^\pm_\delta(\lambda t)b(t)$, then
$$R^\delta_\lambda f=  \rho_\lambda(\lambda-\mathcal{P})f +
 \rho_\lambda(\lambda+\mathcal{P})f.$$ Using
\eqref{3.5}, one finds that for fixed $\delta$ one has the uniform
bounds
$$|\rho_\lambda(\tau)|\le C_N \lambda^{-\delta}(1+|\tau|)^{-N},$$
for every $N$.  Hence, \eqref{3.3} and H\"older's inequality imply
that for every $\lambda\ge 1$
$$\|R^\delta_\lambda f\|_p \le C\|R^\delta_\lambda f\|_2 \le C
\lambda^{\delta(p)-\delta}\|f\|_p,$$ which shows that the
remainder terms $R^\delta_\lambda$ are uniformly bounded when
$\delta>\delta(p)$.

If we let $a(t)=1-b(t)$, so that $a(t)=1$ for $|t|<1$ and $0$ for
$|t|>2$, we would be done if we could prove the same for
$$\tilde S^\delta_\lambda f
=\frac1{2\pi}\int_{-\infty}^\infty \sum_{\pm}\lambda
m_\delta^\pm(\lambda t)a(t) e^{\pm i \lambda t}
 \cos t\mathcal{P}
f\, dt.$$ To do this, we shall make a dyadic decomposition of the
integral.  Fix $\beta\in C^\infty_0(\mathbb{R})$ satisfying
$\beta(t)=0$ $t\notin [1/2,4]$ and $\sum_{-\infty}^\infty
\beta(2^{-j} t)=1$, $t>0$.   We then set
$\beta_0(t)=\sum_{j=0}^\infty \beta(2^{-j}|t|)$ so that $\beta_0$
is smooth and satisfies $\beta_0(t)=0$, $|t|>2$.  We then define
for $j=1,2,\dots$
$$\tilde S^\delta_{\lambda,j} f
=\frac1{2\pi}\int_{-\infty}^\infty \sum_{\pm}\lambda
m_\delta^\pm(\lambda t)\beta(\lambda 2^{-j}|t|)a(t) e^{\pm i
\lambda t}
 \cos t\mathcal{P}
f\, dt,$$ and
$$\tilde S^\delta_{\lambda,0} f
=\frac1{2\pi}\int_{-\infty}^\infty \sum_{\pm}\lambda
m_\delta^\pm(\lambda t)\beta_0(\lambda t)a(t) e^{\pm i \lambda t}
 \cos t\mathcal{P}
f\, dt,$$ so that $\tilde S^\delta_\lambda f= \sum_{j\ge 0}\tilde
S^\delta_{\lambda,j}f$. Note that, because of the support
properties of $a(t)$, $\tilde S^\delta_{\lambda,j}f$ vanishes if
$j$ is larger than a fixed multiple of $\log \lambda$.

We claim that if $\delta>\delta(p)$ is fixed then
\begin{equation}\label{3.6}
\|\tidle S^\delta_{\lambda,j}f\|_p \le C
2^{-(\delta-\delta(p))j}\|f\|_p,\end{equation} where $C$ is
independent of $\lambda$ and $j$.  This would of course complete
the missing step of obtaining the uniform boundedness of the
$\tilde S^\delta_\lambda$ for $\delta>\delta(p)$.

To prove this estimate we shall exploit the fact that the finite
propagation speed of the wave equation mentioned before implies
that the kernels of the operators, $\tilde
S^\delta_{\lambda,j}(x,y)$ must satisfy
$$\tilde S^\delta_{\lambda,j}(x,y)=0, \quad \text{if } \,
\text{dist}(x,y)\ge 8(2^j\lambda^{-1}),$$ since $\cos
t\mathcal{P}$ will have a kernel that vanishes on this set when
$t$ belongs to the support of the integral defining $\tilde
S^\delta_{\lambda,j}$.  Because of this, in order to prove
\eqref{3.6}, it suffices to show that for all geodesic balls
$B_{R_{\lambda,j}}$ of radius $R_{\lambda,j}=\lambda^{-1}2^j$ on
has the bounds
\begin{equation}\label{3.7}
\|\tidle S^\delta_{\lambda,j}f\|_{L^p(B_{R_{\lambda,j}})} \le C
2^{-(\delta-\delta(p))j}\|f\|_p,
\end{equation}
for the $L^p$ norms over $B_{R_{\lambda,j}}$, with $C$, as before,
being independent of $\lambda$ and $j$.  However, by H\"older's
inequality,
$$\|\tidle S^\delta_{\lambda,j}f\|_{L^p(B_{R_{\lambda,j}})} \le C
(\lambda^{-1}2^j)^{\frac{n}p-\frac{n}2} \|\tidle
S^\delta_{\lambda,j}f\|_{L^2(M)},
$$
and so we would be done if we could show that
\begin{equation}\label{3.8}
\|S^\delta_{\lambda,j}f\|_2 \le
C(\lambda^{-1}2^j)^{-\frac{n}p+\frac{n}2}2^{-(\delta-\delta(p))j}\|f\|_p
=\lambda^{\delta(p)}(\lambda 2^{-j})^{\frac12}2^{-\delta j} \,
\|f\|_p .
\end{equation}

To prove this for $j=1,2,\dots$ we note that, by \eqref{3.5}, the
inverse Fourier transform  of $t\to \frac12\lambda
m_\delta(\lambda t)\beta(\lambda 2^{-j}|t|)$ behaves like that of
$\lambda^{-\delta}|t|^{-1-\delta}\beta(\lambda 2^{-j}|t|)$.  Since
the dyadic cutoff localizes to $|t|\approx \lambda^{-1}2^j$, we
conclude that we can write
$$\tilde S^\delta_{\lambda,j}f=
2^{-j\delta}
\rho_{\lambda,j}(\lambda^{-1}2^j(\lambda-\mathcal{P})) +
 2^{-j\delta}
\rho_{\lambda,j}(\lambda^{-1}2^j(\lambda+\mathcal{P})),$$  where
the $\rho_{\lambda,j}$ satisfy the uniform bounds
$$|\rho_{\lambda,j}(\tau)|\le C_N(1+|\tau|)^{-N},$$
for every $N$. Because of this, the estimates \eqref{3.8} with
$j=1,2,\dots$ just follow from \eqref{3.3} with $2^k$ being
replaced by $\lambda 2^{-j}$.  Since the estimate for $j=0$
follows from the same argument, the proof is complete. \qed

\end{document}